\documentclass[12pt]{article} 
\def\R{\hbox{{\rm I}\kern-0.2em{\rm R}\kern0.2em}}%mathematical R for reals
\def\D{\hbox{{\rm I}\kern-0.2em{\rm D}\kern0.2em}}

\def\be{\begin{equation}}
\def\ee{\end{equation}}

\def\({\left(}
\def\){\right)}
\def\[{\left[}
\def\]{\right]}
\def\bc{\begin{center}}
\def\ec{\end{center}}

\begin{document}
{\large \bf Linearization of two dimensional complex-linearizable
systems of second order ordinary differential equations}
\bc{\textbf{Hina M. Dutt}\\
School of Natural Sciences, National University of Sciences and Technology, Campus H-12, 44000, Islamabad, Pakistan\\
hinadutt@yahoo.com}\ec
\bc{\textbf{M. Safdar}\\
School of Mechanical and Manufacturing Engineering, National
University of Sciences and Technology, Campus H-12, 44000,
Islamabad, Pakistan\\
safdar.camp@gmail.com}\ec

{\bf Abstract}. Complex-linearization of a class of systems of
second order ordinary differential equations (ODEs) has already been
studied with complex symmetry analysis. Linearization of this class
has been achieved earlier by complex method, however, linearization
criteria and the most general linearizable form of such systems have
not been derived yet. In this paper, it is shown that the general
\emph{linearizable form of the complex-linearizable systems} of two
second order ODEs is (at most) quadratically semi-linear in the
first order derivatives of the dependent variables. Further,
linearization conditions are derived in terms of coefficients of
system and their derivatives. These linearizable $2$-dimensional
complex-linearizable systems of second order ODEs are characterized
here, by adopting both the real and complex procedures.

\section{Introduction}
Most of the algorithms constructed to solve differential equations
(DEs) with symmetry analysis involve an invertible change of the
dependent and/or independent (point transformations) variables. For
solving nonlinear DEs, symmetry analysis uses a tool called
linearization, which maps them to linear equations under invertible
change of the variables. Linearization procedure requires; the most
general forms of the DEs that could be candidates of linearization
and linearization criteria that ensure existence of invertible
transformations from nonlinear to linear equations. Though
construction of point transformations and finally getting to an
analytic solution of the concerned problem are also involved in
linearization process, these issues are of secondary nature as one
needs to first investigate linearizability of DEs. An explicit
linearizable form and linearization criteria for the scalar second
order ODEs have been derived by Sophus Lie (see, e.g.,
\cite{ibrbook2}). Similarly, linearization of higher order scalar
ODEs and systems of these equations attracted a great deal of
interest and studied comprehensively over the last decade (see,
e.g., \cite{ibr12}-\cite{wafomahomedquad}).

Complex symmetry analysis has been employed to solve certain classes
of systems of nonlinear ODEs and linear PDEs. Of particular interest
here, is linearization of systems of second order ODEs (see, e.g.,
\cite{saj}-\cite{saf2}) that is achieved by complex methods. These
classes are obtained from linearizable scalar and systems of ODEs by
considering their dependent variables as complex functions of a real
independent variable, which when split into the real and imaginary
parts give two dependent variables. In this way, a scalar ODE
produces a system of two coupled equations, with Cauchy-Riemann (CR)
structure on both the equations. These CR-equations appear as
constraint equations that restrict the emerging systems of ODEs to
special subclasses of the general class of such systems. These
subclasses of $2$-dimensional systems of second order ODEs may
trivially be studied with real symmetry analysis, however, they
appear to be nontrivial when viewed from complex approach.
Complex-linearizable (c-linearizable) classes explored earlier
\cite{saj}-\cite{saf2} and studied in this paper provide us means to
extend linearization procedure to $m$-dimensional systems ($m\geq
3$), of $n^{th}$ order ($n\geq 2$) ODEs. Though these classes are
subcases of the general $m$-dimensional systems of $n^{th}$ order
ODEs, their linearization has not been achieved yet, with real
symmetry analysis. Presently symmetry classification and solvability
of higher dimensional systems of higher order ODEs seems to be
exploitable \emph{only with complex symmetry analysis}.

When linearizable scalar second order ODEs are considered complex by
taking the dependent variable as a complex function of a real
independent variable, they lead to c-linearization. The associated
linearization criteria that consist of two equations (see, e.g.,
\cite{ibrbook2}) involving coefficients of the second order
equations and their partial derivatives of (at most) order two, also
yield four constraint equations for the corresponding system of two
ODEs on splitting the complex functions involved, into the real and
imaginary parts. These four equations constitute the c-linearization
criteria \cite{saj}, for the corresponding class of systems of two
second order ODEs. The reason for calling them c-linearization
instead of linearization criteria is that, in earlier works,
explicit \emph{Lie procedure} to obtain linearization conditions of
this class of systems, was not performed after incorporating complex
symmetry approach on scalar ODEs. The \emph{most general form of the
c-linearizable $2$-dimensional linearizable systems} of second order
ODEs is obtained here by real and complex methods. This derivation
shows that the general linearizable forms (obtained by real and
complex procedures) of $2$-dimensional c-linearizable systems of
second order ODEs are identical. Moreover, associated linearization
criteria have been derived, again by adopting both the real and
complex symmetry methods. These linearization conditions are also
shown to be similar whether derived from real Lie procedure
developed for systems or by employing complex symmetry analysis on
scalar ODE. The core result obtained here is refinement of the
c-linearization conditions to linearization criteria for
$2$-dimensional systems of second order ODEs, obtainable from
linearizable complex scalar second order ODEs.

The plan of the paper is as follows. The second section presents
derivation of the linearizable form for the scalar second order ODEs
and Lie procedure to obtain associated linearization criteria. The
subsequent section is on the linearization of $2$-dimensional
c-lineariable systems of second order ODEs, by real and complex
symmetry methods. The fourth section contains some illustrative
examples. The last section concludes the paper.

\section{A subclass of linearizable scalar second order ODEs}
The following point transformations
\begin{eqnarray}
\tilde{x}=\phi(x,u),~~\tilde{u}=\psi(x,u),\label{gentrans}
\end{eqnarray}
where $\phi$ and $\psi$ are arbitrary functions of $x$ and $u$,
yield the most general form of linearizable scalar second order ODEs
\begin{eqnarray}
u^{\prime\prime}+\alpha(x,u)u^{\prime 3}+\beta(x,u)u^{\prime
2}+\gamma(x,u)u^{\prime}+\delta(x,u)=0,\label{cubiclineq}
\end{eqnarray}
with four arbitrary coefficients, that is cubically semi-linear in
the first order derivative of the dependent variable, for derivation
see \cite{ibrbook2}. Restricting these transformations to
\begin{eqnarray}
\tilde{x}=\phi(x),~~\tilde{u}=\psi(x,u),\label{restrans}
\end{eqnarray}
i.e., assuming $\phi_{u}=0$, leads to a quadratically semi-linear
scalar second order ODE that is derived here explicitly. Under
transformations (\ref{restrans}) the first and second order
derivatives of $\tilde{u}(\tilde{x})$ with respect to $\tilde{x}$
read as
\begin{eqnarray}
\tilde{u}^{\prime}=\frac{D \psi(x,u)}{D
\phi(x)}=\lambda(x,u,u^{\prime}),
\end{eqnarray}
and
\begin{eqnarray}
\tilde{u}^{\prime\prime}=\frac{D \lambda(x,u,u^{\prime})}{D
\phi(x)}=\mu(x,u,u^{\prime},u^{\prime\prime}),
\end{eqnarray}
respectively. Here
\begin{eqnarray}
D=\frac{\partial}{\partial x}+u^{\prime}\frac{\partial}{\partial
u}+u^{\prime\prime}\frac{\partial}{\partial u^{\prime}}+\cdots,
\end{eqnarray}
is the total derivative operator. Inserting the total derivative
operator in both the above equations leads us to the following
\begin{eqnarray}
\tilde{u}^{\prime}=\frac{\psi_{x}+u^{\prime}\psi_{u}}{\phi_{x}},
\end{eqnarray}
and
\begin{eqnarray}
\tilde{u}^{\prime\prime}=\frac{\phi_{x}(\psi_{xx}+2u^{\prime}\psi_{xu}+u^{\prime
2}\psi_{uu}+u^{\prime\prime}\psi_{u})-\phi_{xx}(\psi_{x}+u^{\prime}\psi_{u})}{\phi^{3}_{x}},\label{transeq1}
\end{eqnarray}
respectively. Equating (\ref{transeq1}) to zero, i.e., considering
$\tilde{u}^{\prime\prime}=0,$ leaves a quadratically semi-linear ODE
of the form
\begin{eqnarray}
u^{\prime\prime}+a(x,u)u^{\prime
2}+b(x,u)u^{\prime}+c(x,u)=0,\label{quadlineq}
\end{eqnarray}
with the coefficients
\begin{eqnarray}
a(x,u)=\frac{\psi_{uu}}{\psi_{u}},~~b(x,u)=\frac{2\phi_{x}\psi_{xu}-\psi_{u}\phi_{xx}}{\phi_{x}\psi_{u}},~~c(x,u)=\frac{\phi_{x}\psi_{xx}-\psi_{x}\phi_{xx}}{\phi_{x}\psi_{u}}.\label{coeffquadsca}
\end{eqnarray}
The quadratic nonlinear (in the first derivative) equation
(\ref{quadlineq}) with three coefficients (\ref{coeffquadsca}) is a
subcase of the general linearizable (cubically semi-linear) second
order ODE (\ref{cubiclineq}).

Now for the derivation of Lie linearization criteria of nonlinear
equation (\ref{quadlineq}), we start with a re-arrangement
\begin{eqnarray}
\psi_{uu}=a(x,u)\psi_u~,~~~~~~~~~~~~~~~~\nonumber\\
2\psi_{xu}=\phi_x^{-1}\psi_u\phi_{xx}+b(x,u)\psi_u~,\nonumber\\
\psi_{xx}=\phi_x^{-1}\psi_x\phi_{xx}+c(x,u)\psi_u~.
\end{eqnarray}
of the relations (\ref{coeffquadsca}). Equating the mixed
derivatives of $\psi$, such that $(\psi_{xu})_u=(\psi_{uu})_x$ and
$(\psi_{xu})_x=(\psi_{xx})_u,$ we find
\begin{eqnarray}
b_u-2a_x=0,\label{lincriquadsca}
\end{eqnarray}
and
\begin{eqnarray}
\phi_x^{-2}(2\phi_x\phi_{xx}-3\phi_{xx}^2)=4(c_u+ac)-(2b_x+b^2).\label{lincriquadsca1}
\end{eqnarray}
As $\phi_u=0$, differentiating (\ref{lincriquadsca1}) with respect
to $u$, simplifies it to
\begin{eqnarray}
c_{uu}-a_{xx}-a_xb+a_uc+c_ua=0.\label{lincriquadsca2}
\end{eqnarray}
Equations (\ref{lincriquadsca}) and (\ref{lincriquadsca2})
constitute the linearization criteria for the scalar second order
quadratically semi-linear ODEs \cite{sim}.

\section{Linearizable two dimensional c-linearizable systems of second order ODEs}
To obtain the most general form of a \emph{linearizable
$2$-dimensional c-linearizable system} of second order ODEs we
consider two approaches viz, by complex symmetry analysis and by
real symmetry analysis.
\subsection{Use of complex symmetry analysis}
Suppose $u(x)$ in (\ref{quadlineq}) be complex function of a real
variable $x$ i.e., $u(x)=y(x)+iz(x)$. Further assume that
\begin{eqnarray}
a(x,u)=a_1(x,y,z)+ia_2(x,y,z)~,\nonumber\\
b(x,u)=b_1(x,y,z)+ib_2(x,y,z)~,\nonumber\\
c(x,u)=c_1(x,y,z)+ic_2(x,y,z)~.\label{splitcoeff}
\end{eqnarray}
This converts the scalar ODE (\ref{quadlineq}) to a system of two second order ODEs of the form
\begin{eqnarray}
y''+a_1y'^2-2a_2y'z'-a_1z'^2+b_1y'-b_2z'+c_1=0~,\nonumber\\
z''+a_2y'^2+2a_1y'z'-a_2z'^2+b_2y'+b_2z'+c_2=0~,\label{splitsys}
\end{eqnarray}
with the coefficients $a_j, b_j, c_j;(j=1,2),$ satisfying the
CR-equations
\begin{eqnarray}
a_{1,y}=a_{2,z},~~a_{1,z}=-a_{2,y},\nonumber\\
b_{1,y}=b_{2,z},~~b_{1,z}=-b_{2,y},\nonumber\\
c_{1,y}=c_{2,z},~~c_{1,z}=-c_{2,y}.\label{comcoeffs54}
\end{eqnarray}
Moreover, conditions (\ref{lincriquadsca}) and
(\ref{lincriquadsca2}) can now be converted into a set of four
equations
\begin{eqnarray}
2a_{1,x}-b_{1,y}=0,~~~~~~~~~~~~~~~~~~~~~~~~~~~~~~~~~~~~~~~~~~~~~~~~~\label{c1}\\
2a_{2,x}+b_{1,z}=0,~~~~~~~~~~~~~~~~~~~~~~~~~~~~~~~~~~~~~~~~~~~~~~~~~\label{c2}\\
c_{1,zz}+a_{1,xx}+a_{1,x}b_1-a_{2,x}b_2-(a_2c_1)_{,z}-(a_1c_2)_{,z}=0,\label{c3}\\
c_{2,yy}-a_{2,xx}-a_{2,x}b_1-a_{1,x}b_2+(a_2c_1)_{,y}+(a_1c_2)_{,y}=0,\label{c4}
\end{eqnarray}
by splitting the complex coefficients (\ref{comcoeffs54}) into the
real and imaginary parts.

As evident from \cite{saj}, such a (complex) procedure leads us to
c-linearization of systems of ODEs. Our claim here is that equations
(\ref{c1}-\ref{c4}) are actually the linearization conditions
despite of being just the c-linearization conditions for system
(\ref{splitsys}). In order to prove this fact, we now adopt real
symmetry analysis in the next subsection to derive the lineariation
conditions for system (\ref{splitsys}).

\subsection{Use of real symmetry analysis}
The previous work on c-linearizable \cite{saj,saf2} and their
linearizable subclass of systems \cite{saf1,saf4} of second order
ODEs reveals that point transformations of the form
\begin{eqnarray}
\tilde{x}=\phi(x),~~\tilde{y}=\psi_{1}(x,y,z),~~\tilde{z}=\psi_{2}(x,y,z),\label{resfibrpre}
\end{eqnarray}
where
\begin{eqnarray}
\psi_{1,y}=\psi_{2,z},~~\psi_{2,y}=-\psi_{1,z},\label{CREs0}
\end{eqnarray}
i.e., $\psi_{j}$, for $j=1,2,$ satisfy the CR-equations that involve
derivatives with respect to both the dependent variables, linearizes
the c-linearizable systems. Notice that (\ref{resfibrpre}) are
obtainable from (\ref{restrans}) that is a subclass of
(\ref{gentrans}). These transformations map the first and second
order derivatives as
\begin{eqnarray}
\tilde{y}^{\prime}=\frac{D\psi_{1}}{D\phi}=\lambda_{1}(x,y,z,y^{\prime},z^{\prime}),~~\tilde{z}^{\prime}=\frac{D\psi_{2}}{D\phi}=\lambda_{2}(x,y,z,y^{\prime},z^{\prime}),
\end{eqnarray}
and
\begin{eqnarray}
\tilde{y}^{\prime\prime}=\frac{D\lambda_{1}}{D\phi}=\mu_{1}(x,y,z,y^{\prime},z^{\prime},y^{\prime\prime},z^{\prime\prime}),~~\tilde{z}^{\prime\prime}=\frac{D\lambda_{2}}{D\phi}=\mu_{2}(x,y,z,y^{\prime},z^{\prime},y^{\prime\prime},z^{\prime\prime}),
\end{eqnarray}
where
\begin{eqnarray}
D=\frac{\partial}{\partial x}+y^{\prime}\frac{\partial}{\partial
y}+z^{\prime}\frac{\partial}{\partial
z}+y^{\prime\prime}\frac{\partial}{\partial
y^{\prime}}+z^{\prime\prime}\frac{\partial}{\partial
z^{\prime}}+\cdots.
\end{eqnarray}
Inserting the total derivative operator in the above equations and
simplifying, we arrive at the following $2$-dimensional system
\begin{eqnarray}
y^{\prime\prime}+\alpha_{1}y^{\prime
2}-2\alpha_{2}y^{\prime}z^{\prime}+\alpha_{3}z^{\prime
2}+\beta_{1}y^{\prime}-\beta_{2}z^{\prime}+\gamma_{1}=0,\nonumber\\
z^{\prime\prime}+\alpha_{4}y^{\prime
2}+2\alpha_{5}y^{\prime}z^{\prime}+\alpha_{6}z^{\prime
2}+\beta_{3}y^{\prime}+\beta_{4}z^{\prime}+\gamma_{2}=0,\label{linform0}
\end{eqnarray}
where
\begin{eqnarray*}
\alpha_{1}=\phi_x\Delta^{-1}(\psi_{2,z}\psi_{1,yy}-\psi_{1,z}\psi_{2,yy}),~~\alpha_{2}=\phi_x\Delta^{-1}(\psi_{1,z}\psi_{2,yz}-\psi_{2,z}\psi_{1,yz}),\nonumber\\
\alpha_{3}=\phi_x\Delta^{-1}(\psi_{2,z}\psi_{1,zz}-\psi_{1,z}\psi_{2,zz}),~~\alpha_{4}=\phi_x\Delta^{-1}(\psi_{1,y}\psi_{2,yy}-\psi_{2,y}\psi_{1,yy}),\nonumber\\
\alpha_{5}=\phi_x\Delta^{-1}(\psi_{1,y}\psi_{2,yz}-\psi_{2,y}\psi_{1,yz}),~~\alpha_{6}=\phi_x\Delta^{-1}(\psi_{1,y}\psi_{2,zz}-\psi_{2,y}\psi_{1,zz}),
\end{eqnarray*}
\begin{eqnarray*}
\beta_{1}=2\phi_x\Delta^{-1}(\psi_{2,z}\psi_{1,xy}-\psi_{1,z}\psi_{2,xy})-\frac{\phi_{xx}}{\phi_x},~~\beta_{2}=2\phi_x\Delta^{-1}(\psi_{1,z}\psi_{2,xz}-\psi_{2,z}\psi_{1,xz}),\nonumber\\
\beta_{3}=2\phi_x\Delta^{-1}(\psi_{1,y}\psi_{2,xy}-\psi_{2,y}\psi_{1,xy}),~~\beta_{4}=2\phi_x\Delta^{-1}(\psi_{1,y}\psi_{2,xz}-\psi_{2,y}\psi_{1,xz})-\frac{\phi_{xx}}{\phi_x},\nonumber\\
\end{eqnarray*}
and
\begin{eqnarray}
\gamma_1=\Delta^{-1}(\phi_{x}\psi_{1,y}\psi_{1,xx}-\psi_{1,x}\psi_{1,y}\phi_{xx}-\phi_x\psi_{1,z}\psi_{2,xx}+\psi_{1,z}\psi_{2,x}\phi_{xx})~,\nonumber\\
\gamma_2=\Delta^{-1}(\phi_{x}\psi_{1,z}\psi_{1,xx}-\psi_{1,x}\psi_{1,z}\phi_{xx}+\phi_x\psi_{1,y}\psi_{2,xx}+\psi_{1,y}\psi_{2,x}\phi_{xx})~,\label{coeffsquadsys}
\end{eqnarray}
where
\begin{eqnarray}
\Delta=\phi_{x}(\psi_{1,y}\psi_{2,z}-\psi_{1,z}\psi_{2,y})\neq0~,
\end{eqnarray}
is the Jacobian of the transformation (\ref{resfibrpre}). The coefficients (\ref{coeffquadsca}) of the scalar ODE
(\ref{quadlineq}) split into the coefficients of the corresponding
$2$-dimensional system of second order ODEs. This happens due to
presence of the complex dependent function $u$, in the coefficients
(\ref{coeffquadsca}). The restricted fibre preserving
transformations (\ref{resfibrpre}) used to derive the linearizable
form (\ref{linform0}), are obtainable from the complex
transformations (\ref{restrans}) that are employed to deduce
(\ref{quadlineq}). Therefore, transformations (\ref{resfibrpre})
along with (\ref{CREs0}) appear to be the real and imaginary parts
of complex transformation (\ref{restrans}), they reveal the
correspondence of the linearizable forms of $2$-dimensional systems
and scalar complex ODEs. The CR-equations are not yet incorporated
in the linearizable form (\ref{linform0}). Insertion of the
CR-equations (\ref{CREs0}) and their derivatives
\begin{eqnarray}
\psi_{1,yy}=\psi_{2,yz}=-\psi_{1,zz},\nonumber\\
\psi_{2,zz}=\psi_{1,yz}=-\psi_{2,yy},\label{CREs1}
\end{eqnarray}
brings out the correspondence between the coefficients
(\ref{coeffquadsca}) of the complex linearizable ODEs
(\ref{quadlineq}) and coefficients (\ref{coeffsquadsys}) of the
system (\ref{linform0}). Employing (\ref{CREs0}) and (\ref{CREs1})
the coefficients (\ref{coeffsquadsys}) reduces to \emph{only six
arbitrary} coefficients that read as
\begin{eqnarray}
\alpha_{1}=-\alpha_{3}=\alpha_{5}=a_{1},~~\alpha_{2}=\alpha_{4}=-\alpha_{6}=a_{2},\nonumber\\
\beta_{1}=\beta_{4}=b_{1},~~\beta_{2}=\beta_{3}=b_{2},~~\gamma_{1}=c_{1},~~\gamma_{2}=c_{2}.\label{samecoeff}
\end{eqnarray}
Here the coefficients $a_{j},~b_{j}$ and $c_{j}$ are the real and
imaginary parts of the complex coefficients (\ref{coeffquadsca}).
The linearizable form of systems derived in this section by real
method appears to be the same as one obtains by splitting the
corresponding form of the scalar complex equation (\ref{quadlineq}).
This analysis leads us to the following theorem.

\textbf{Theorem 1.} \emph{The most general form of the linearizable
two dimensional c-linearizable systems of second order ODEs is
quadratically semi-linear}.

\subsubsection{Sufficient conditions for the linearization of a c-linearizable system}
Consider the most general form of the c-linearizable $2$-dimensional
systems of second order ODEs (\ref{splitsys}), with constraint
equations (\ref{comcoeffs54}). Rewriting the coefficients of the
system (\ref{splitsys}) in the form
\begin{eqnarray}
a_1=\Delta^{-1}\phi_x(\psi_{1,y}\psi_{1,yy}+\psi_{1,z}\psi_{1,yz}),~~~~~~~~~~~~~~~~~~~~~~~~~~~~~~~~~~~~~~~\nonumber\\
a_2=\Delta^{-1}\phi_x(\psi_{1,z}\psi_{1,yy}+\psi_{1,y}\psi_{1,yz}),~~~~~~~~~~~~~~~~~~~~~~~~~~~~~~~~~~~~~~~\nonumber\\
b_1=2\Delta^{-1}\phi_x(\psi_{1,y}\psi_{1,xy}+\psi_{1,z}\psi_{1,xz})-{\phi_{xx}\over \phi_x},~~~~~~~~~~~~~~~~~~~~~~~~~~~~~\nonumber\\
b_2=2\Delta^{-1}\phi_x(\psi_{1,z}\psi_{1,xy}+\psi_{1,y}\psi_{1,xz}),~~~~~~~~~~~~~~~~~~~~~~~~~~~~~~~~~~~~~~\nonumber\\
c_1=\Delta^{-1}(\phi_{x}\psi_{1,y}\psi_{1,xx}-\psi_{1,x}\psi_{1,y}\phi_{xx}-\phi_x\psi_{1,z}\psi_{2,xx}+\psi_{1,z}\psi_{2,x}\phi_{xx}),\nonumber\\
c_2=\Delta^{-1}(\phi_{x}\psi_{1,z}\psi_{1,xx}-\psi_{1,x}\psi_{1,z}\phi_{xx}+\phi_x\psi_{1,y}\psi_{2,xx}+\psi_{1,y}\psi_{2,x}\phi_{xx}).\label{coeff2}
\end{eqnarray}
For obtaining the sufficient linearizability conditions of (\ref{splitsys}), we have
to solve compatibility problem, that has already been solved for the
scalar equations earlier in this work, for the set of equations
(\ref{coeff2}). It is an over determined system of partial
differential equations for the functions $\phi, \psi_1$ and $\psi_2$
with known $a_j, b_j, c_j$.\\
The system (\ref{coeff2}) gives us
\begin{eqnarray}
\psi_{1,yy}&=&\psi_{1,y}a_1+\psi_{1,z}a_2~,\nonumber\\
\psi_{1,yz}&=&\psi_{1,z}a_1-\psi_{1,y}a_2~,\nonumber\\
\psi_{1,xy}&=&{1\over 2}(\psi_{1,y}b_1+\psi_{1,z}b_2+\psi_{1,y}{\phi_{xx}\over \phi_x})~,\nonumber\\
\psi_{1,xz}&=&{1\over
2}(\psi_{1,z}b_1-\psi_{1,y}b_2+\psi_{1,z}{\phi_{xx}\over
\phi_x})~,\nonumber\\
\psi_{1,xx}&=&\psi_{1,y}c_{1}+\psi_{1,z}c_2+\psi_{1,x}{\phi_{xx}\over \phi_x}~,\nonumber\\
\psi_{2,xx}&=&\psi_{1,y}c_{2}-\psi_{1,z}c_1+\psi_{2,x}{\phi_{xx}\over
\phi_x}~.\nonumber
\end{eqnarray}

The compatibility of the system (\ref{coeff2}) first requires to compute partial
derivatives
\begin{eqnarray}
\Delta_x&=&2\Delta{\phi_{xx}\over \phi_x}+\Delta b_1~,\nonumber\\
\Delta_y&=&2\Delta a_1~,\nonumber\\
\Delta_z&=&-2\Delta a_2~,\nonumber
\end{eqnarray}
of the Jacobian. Comparing the mixed derivatives
$(\Delta_y)_z=(\Delta_z)_y$, $(\Delta_x)_y=(\Delta_y)_x$ and
$(\Delta_x)_z=(\Delta_z)_x$, we obtain
\begin{eqnarray}
a_{1,z}+a_{2,y}=0~,\label{constraint1}\\
2a_{1,x}-b_{1,y}=0~,\label{constraint2}\\
2a_{2,x}+b_{1,z}=0~,\label{constraint3}
\end{eqnarray}
respectively.
Equating the mixed derivatives $(\psi_{1,yy})_z=(\psi_{1,yz})_y$, $(\psi_{1,yy})_x=(\psi_{1,xy})_y$, $(\psi_{1,xx})_y=(\psi_{1,xy})_x$ , $(\psi_{1,xx})_z=(\psi_{1,xz})_x$, $(\psi_{1,xy})_z=(\psi_{1,xz})_y$, $(\psi_{2,xx})_y=(\psi_{2,xy})_y$ and $(\psi_{2,xx})_z=(\psi_{2,xz})_x$ gives us
\begin{eqnarray}
a_{1,y}-a_{2,z}=0,~~~~~~~~~~~~~~~~~~~~~~~~~~~~~~~~~~~~~~~~~~~~~~~~~~~\label{constraint4}\\
b_{2,y}+b_{1,z}=0,~~~~~~~~~~~~~~~~~~~~~~~~~~~~~~~~~~~~~~~~~~~~~~~~~~~\label{constraint6}\\
b_{2,z}-b_{1,y}=0,~~~~~~~~~~~~~~~~~~~~~~~~~~~~~~~~~~~~~~~~~~~~~~~~~~~\label{constraint7}\\
c_{2,z}-c_{1,y}=0,~~~~~~~~~~~~~~~~~~~~~~~~~~~~~~~~~~~~~~~~~~~~~~~~~~~\label{constraint8}\\
c_{2,y}+c_{1,z}=0,~~~~~~~~~~~~~~~~~~~~~~~~~~~~~~~~~~~~~~~~~~~~~~~~~~~\label{constraint9}\\
c_{1,zz}+a_{1,xx}+a_{1,x}b_1-a_{2,x}b_2-(a_2c_1)_{,z}-(a_1c_2)_{,z}=0,\label{constraint11}\\
c_{2,yy}-a_{2,xx}-a_{2,x}b_1-a_{1,x}b_2+(a_1c_2)_{,y}-(a_2c_1)_{,y}=0.\label{constraint12}\\
\end{eqnarray}
Note that $(\psi_{1,yz})_x-(\psi_{1,xz,})_y=0$ and $(\psi_{1,xy})_z-(\psi_{1,yz})_x=0$ are satisfied.
Also (\ref{constraint1}), (\ref{constraint4})-(\ref{constraint9}) are CR-equations for the
coefficients $a_j,$ $b_j,$ $c_j$. Therefore, the solution of the
compatibility problem of the system (\ref{coeff2}), provides
CR-constraints on the coefficients of (\ref{splitsys}) and the
linearization conditions.

\textbf{Theorem 2.} \emph{A two dimensional c-linearizable system of
second order ODEs of the form} (\ref{splitsys}) \emph{is
linearizable if and only if its coefficients satisfy the
CR-equations and conditions} (\ref{constraint2}),
(\ref{constraint3}), (\ref{constraint11}), (\ref{constraint12}).

These are the same conditions that are already obtained
(\ref{c1}-\ref{c4}), by employing complex analysis, i.e., splitting
the linearization conditions associated with the base scalar
equation (\ref{quadlineq}), into the real and imaginary parts.

\textbf{Corollary.} \emph{The c-linearization conditions for a two
dimensional system of quadratically semi-linear second order ODEs
are the linearization conditions.}

\section{Examples}
We present some examples to illustrate our results.\\\\
\textbf{1.}~The $2$-dimensional system of second order ODEs
\begin{eqnarray}
y''-({2y\over y^2+z^2})y'^2-2({2z\over y^2+z^2})y'z'+({2y\over y^2+z^2})z'^2-{2\over x}y'-{2y\over x^2}=0~,\nonumber\\
z''+({2z\over y^2+z^2})y'^2-2({2z\over y^2+z^2})y'z'-({2z\over y^2+z^2})z'^2-{2\over x}z'-{2z\over x^2}=0~.\label{example1}
\end{eqnarray}
is of the same form as (\ref{splitsys}) with
\begin{eqnarray}
a_1={-2y\over y^2+z^2}~,\quad a_2={2z\over y^2+z^2}~,\quad
b_1={-2\over x}~,\quad b_2=0~,\quad c_1={-2y\over x^2}~,\quad
c_2={-2z\over x^2}~.\label{coeffexample1}
\end{eqnarray}
One can easily verify that (\ref{coeffexample1}) satisfy the conditions (\ref{constraint2}), (\ref{constraint3}), (\ref{constraint11}), (\ref{constraint12}) and CR-equations w.r.t $y$ and $z$. So the system of ODEs (\ref{example1}) is linearizable. The transformation
\begin{eqnarray}
t=x~,\quad u={y\over x(y^2+z^2)}~,\quad v={-z\over x(y^2+z^2)}~,
\end{eqnarray}
reduces the nonlinear system (\ref{example1}) to the linear system $u''=0~,v''=0~.$
\\\\
\textbf{2.} Consider the following system of nonlinear ODEs
\begin{eqnarray}
y''-{1\over f(y,z)}(y'^2\cos y\sin y-z'^2\cos y\sin y-2y'z'\cosh z\sinh z)+{2y'\over x}=0~,\nonumber\\
z''-{1\over f(y,z)}(y'^2\cosh z\sinh z-z'^2\cosh z\sinh z+2y'z'\cosh y\sinh y)+{2z'\over x}=0~,\label{example2}
\end{eqnarray}
where $f(y,z)=\sin^2y\cosh^2z+\cos^2y\sinh^2y$, and the
coefficients satisfy the CR-constraint and linearization conditions (\ref{constraint2}), (\ref{constraint3}), (\ref{constraint11}), (\ref{constraint12}). Hence Theorem 2
guarantees that system (\ref{example2}) can be transformed to system
of linear equations $u''=0~,~v''=0$. The linearizing transformations
in this case are
\begin{eqnarray}
t=x~,\quad u=x\cos y\cosh z~,\quad v=-x\sin y\sinh z~.
\end{eqnarray}\\\\
\textbf{3.} Consider the anisotropic oscillator system
\begin{eqnarray}
y''+f(x)y=0~,\nonumber\\
z''+g(x)z=0~.\label{example3}
\end{eqnarray}
In \cite{mahomedqadircubic} it is shown that system (\ref{example3})
is reducible to the free particle system ($u''=0~,~v''=0$) provided
$f=g$. Our c-linearization criteria also leads to the same
condition, i.e. $f=g$.

\section{Conclusion}
C-linearization of $2$-dimensional systems of second order ODEs is
achieved earlier by considering the scalar second order linearizable
ODEs as complex. Their associated linearization criteria are
separated into the real and imaginary parts due to complex functions
involved. In this work, the c-linearization and linearization are
shown to be two different criteria for a $2$-dimensional systems of
second order ODEs. Linearizable form of such c-linearizable systems
has been derived and it is shown to be quadratically semi-linear in
the first order derivatives. Moreover, complex linearization
criteria have been refined to linearization criteria for such
$2$-dimensional systems that are linearizable due to their
correspondence with the complex scalar ODEs.

Earlier in this work, c-linearizable classes of systems of ODEs are
claimed to be non-trivial, when viewed from complex approach. The
reason for calling them non-trivial is that the concept of
c-linearization of systems of ODEs is extendable to $m$-dimensional
systems of $n^{th}$ order ODEs. The simplest procedure that might
lead us to linearization of $m$-dimensional system of second order
ODEs, is to iteratively complexify a scalar second order
linearizable ODE. Therefore, complex symmetry analysis needs to be
extended to $2$- and $3$-dimensional systems of third and second
order ODEs, respectively, in order to derive the general
linearization results mentioned above. Likewise, complex symmetry
analysis may lead us to algebraic classification of the higher
dimensional systems of higher order ODEs.

\section*{Acknowledgment} The authors are most grateful to Asghar Qadir for useful comments and discussion on this work.

\end{document}